\input amstex 
\documentstyle{amsppt}
\input bull-ppt
\keyedby{bull427/car}

\topmatter
\cvol{29}
\cvolyear{1993}
\cmonth{October}
\cyear{1993}
\cvolno{2}
\cpgs{223-227}
\title A new measure of growth\\ 
for countable-dimensional algebras\endtitle
\author John Hannah and K. C. O'Meara\endauthor
\shorttitle{Growth for Countable-Dimensional Algebras}
\address Department of Mathematics, University of 
Canterbury, Christchurch 1, 
New  Zealand\endaddress
\ml\nofrills{\it E-mail address}:\enspace 
jjh\@math.canterbury.ac.nz
\mlpar
{\it E-mail address}:\enspace 
kco\@math.canterbury.ac.nz\endml
\date October 22, 1992\enddate
\abstract A new dimension function on 
countable-dimensional algebras (over a
field) is described. Its dimension values for finitely 
generated algebras
exactly fill the unit interval $[0,1]$. Since the free 
algebra on two
generators turns out to have dimension 0 (although 
conceivably some Noetherian 
algebras might have positive dimension!), this dimension 
function promises to
distinguish among algebras of infinite 
GK-dimension.\endabstract
\endtopmatter

\document

\heading 1. Introduction\endheading

In recent times, the most prominent dimension used in the 
study of algebras
has been the Gelfand-Kirillov dimension (GK-dimension), 
which measures the
``growth of an algebra in terms of generators'' (see 
\cite2). Here we present
another view of ``growth of an algebra'', based on certain 
infinite matrix 
representations. By an {\it algebra\/} we shall always 
mean an associative 
algebra over a field, with an identity element.

It would be unthinkable that one could have a serious 
study of 
{\it finite-dimen\-sional\/} 
algebras without ever resorting to finite matrix
representations. In the theory of {\it 
infinite-dimensional\/} algebras,
however, {\it infinite\/} matrix representations have 
played only a very minor
role. One reason for this, perhaps, is that the nice 
``arithmetic'' functions
provided by a finite matrix representation---such as 
trace, determinant, rank,
etc.---would appear to have no cousins in the 
infinite-dimensional case.
However, a recent and surprising result by Goodearl, 
Menal, and Moncasi
\cite{1, Proposition 2.1} offers fresh hope for infinite 
matrix representations
of {\it countable-dimensional\/} algebras $A$ over a field 
$F$: it says that
such $A$ can be embedded in the algebra $B(F)$ of all 
$\omega \times \omega$
matrices over $F$ which are simultaneously row-finite and 
column-finite. (Note
$\omega\times \omega$ matrices are just $\aleph_0\times 
\aleph_0$ matrices with
their rows and columns ordered in the standard way.) This 
result has been the
inspiration for our work. For in any such representation 
of $A$, the elements
of $A$ now have {\it all\/} their nonzero entries 
relatively close to the main
diagonal. This raises the question of just how closely 
these nonzero entries
can be squeezed to the main diagonal for a suitable 
embedding. To help
quantify this, we introduce the notion of a {\it growth 
curve\/} for an element
$x\in B(F)$. We say that a function $g\:\bold N\rightarrow 
\bold R^+$ is a
growth curve for $x\in B(F)$ if, for each $n\in \bold N$, 
$x(n,i)=0=x(i,n)$
for all $i>n+g(n)$. In other words, $g(n)$ gives a bound 
on the ``bandwidth''
of $x$ at the $(n,n)$ position, if we interpret bandwidth 
as in Figure
1. (There are other interpretations of ``bandwidth'', of 
course.)

\fg{15pc}
\caption{Figure 1}

We say that $x\in B(F)$ has {\it at most order\/} $g(n)$ 
{\it growth\/} (or that
$x$ has $\roman O(g(n))$ growth) if there is a constant 
$c>0$ such that the function
$cg(n)$ is a growth curve for $x$. If $A$ is a subalgebra 
of $B(F)$ and every
$x\in A$ has $\roman O(g(n))$ growth, then we say that the 
algebra $A$ itself
has $\roman O(g(n))$ growth (but notice that the constant 
$c$ in $cg(n)$ will
depend on the particular $x\in A)$. To say $A$ has {\it 
linear growth\/} means
$A$ has $\roman O(n)$ growth.

\heading 2. Statement of main results\endheading

Our first theorem is an improvement on the Goodearl, 
Menal, and Moncasi
embedding mentioned earlier.
\proclaim{Theorem 1} Every countable-dimensional algebra 
$A$ over a field $F$
can be embedded in $B(F)$ as a subalgebra of linear 
growth.\endproclaim

The proof relies on two key results. The first says that 
for any finite number
$k$ of linear transformations of a countable-dimensional 
vector space, there is
a simultaneous block tridiagonal matrix form of them in 
which the sizes of the
diagonal blocks are
$$1,\ 2k+1,\ (2k+1)^2,\dots, (2k+1)^n,\dots\.$$
The second result is the fact that any 
countable-dimensional algebra can be
embedded in some finitely generated algebra \cite3. We 
stress that all our 
algebra embeddings are required to preserve the identity. 
Full details of the
proof will appear elsewhere.

Can we improve on Theorem 1? The answer is ``no'' in 
general, or at least not
in terms of $\roman O(n^r)$ growth for any $r<1$. But 
there is an interesting
range of ``sublinear'' growths which we can identify for 
individual algebras.
For each real number $r\in[0,1]$ let
$$G(r)=\{x\in B(F)\,|\,x\text{ has }\roman O(n^r)\text{ 
growth}\}\.$$
These $G(r)$ are subalgebras of $B(F)$, and Theorem 1 says 
each countable-
dimensional algebra $A$ can be embedded in $G(1)$. We take 
the ``least'' $r$ for
which $A$ embeds in $G(r)$ as our new dimension for $A$.

\dfn{Definition} The {\it bandwidth dimension\/} of a 
countable-dimensional 
algebra $A$  ({over a field $F)$ is\/}
$$\inf\{r\in \bold R,\ r\geq0\,|\,A\text{ embeds in 
}G(r)\}$$
or, equivalently,
$$\inf\{r\in \bold R,\ r\geq0\,|\,A\text{ embeds in $B(F)$ 
with $\roman O(n^r)$
growth}\}\.$$
\enddfn

Our second theorem completely describes the possible range 
of values of our 
dimension function on finitely generated algebras and, 
hence, also on 
countable-dimensional algebras (by Theorem 1). It suggests 
that bandwidth 
dimension is somewhat ``smoother'' than $GK$-dimension, 
whose corresponding
range of values on finitely generated algebras is
$$0,\ 1,\ 2,\ \text{any real }r>2,\text{ and }+\infty \.$$
\proclaim{Theorem 2} For any field $F$, bandwidth 
dimensions of finitely
generated algebras over $F$ exactly fill the unit interval 
$[0,1]$.
\endproclaim

The proof (over 20 pages) will appear elsewhere in a paper 
by the second
author. A sketch is provided later.

It is not difficult to show that the free algebra 
$F\{x,y\}$ on two generators
embeds in the algebra $G(0)$ of finite bandwidth matrices 
and so has zero
growth. Accordingly, its bandwidth dimension is the {\it 
smallest\/} possible
value, namely, 0. In sharp contrast, from the $GK$ point 
of view this algebra
has exponential growth and, therefore, the {\it largest 
possible 
$GK$-dimension\/},
namely, $+\infty $. This suggests that bandwidth dimension 
is quite different
from $GK$-dimension.

Bandwidth dimension behaves as one would hope with respect 
to subalgebras,
finite subdirect products, and finite matrix algebras. For 
instance, an 
algebra $A$ and its matrix algebra $M_n(A)$ have the same 
dimension.
However, bandwidth dimension behaves poorly on factor 
algebras---but this is
the expected price we pay for not giving the free algebra 
the largest possible
dimension.

There are good signs that certain bandwidth dimension 
values may be reflected
in interesting purely ring-theoretic properties. As an 
illustration of how a
growth curve restriction can result in quite strong ring 
properties, we quote:
\proclaim{Proposition} Suppose $R$ is a \RM(von 
Neumann\RM) regular, right
self-injective ring.
\roster 
\item If $R$ has linear growth over some field $F$ 
\RM(i.e., $R$ embeds in
$G(1))$, then $R$ must be of Type $I_f$.
\item If $R$ has zero growth \RM(i.e., $R$ embeds in 
$G(0))$, then $R$ must have
bounded index of nilpotence.
\endroster
\endproclaim

\heading 3. Sketch of the proof of Theorem 2\endheading

Let $r\in[0,1]$. We shall construct a finitely generated 
algebra $A$ of
bandwidth dimension $r$. The case $r=0$ is trivial, and 
the case $r=1$ turns out
to be taken care of by any ``purely infinite'' $A$, that 
is, $A\cong A\oplus A$
as right $A$-modules. Therefore, we can assume $0<r<1$. 
Let $t=r/(1-r)$, and
for each $k\in \bold N$ let $n_k=[k^t]$ where $[\ ]$ 
denotes the integer part.
In order to absorb some of the flavor of the proof, 
consider first the natural
copy $R$ of $\prod^\infty _{k=1}M_{n_k}(F)$ inside $B(F)$, 
that is, $R$
consists of all the block-diagonal matrices of the form 
shown in Figure 2.
Let 
$J_k$ be the $k$\<th diagonal block for each $k\in \bold 
N$ (so $J_k
\cong M_{n_k}(F))$. The $n_k$ have been designed so that $R
\subseteq G(r)$ but $R\nsubseteq G(s)$ for any $s<r$. Even 
though $R$ is not
countable dimensional, we can still talk of its bandwidth 
dimension (as an
abstract algebra). Although at first sight this looks like 
it ought to be $r$,
in fact the bandwidth dimension is 0 because $R$ can be 
embedded in $G(s)$ for
any $0<s<r$, simply by ``stretching'' out the above 
representation and
repeating blocks often enough. Preventing this stretching 
is one of the key 
problems addressed in the proof. In the case of $R$, it 
turns out
that the algebra $B$ generated by $R$ and the standard 
one-dimensional 
infinite shift is sufficiently ``rigid'' to have bandwidth 
dimension exactly
$r$. But, of course, $B$ is far from being finitely 
generated!


\fg{15pc}
\caption{Figure 2}

In the proof of Theorem 2, we construct $A$ as an 
8-generator subalgebra of
$G(r)$ with the following:
\proclaim{Key Property} $A$ contains each $J_k$, and the 
number of products of
generators needed to obtain each of the standard matrix 
units in $J_k$ grows
essentially logarithmically in $k$\RM; in fact, we can get 
by with $\roman O((
\log k)^2)$ products.\endproclaim

This property is achieved by including various infinite 
shift matrices among
the generators of $A$ and also a block diagonal matrix of 
strategically placed
$2^m\times 2^m$ ``binary search'' matrices for various 
$m$. For elements $a$,
$b$ of a general ring $T$, we say $\gamma\in T$ is a {\it 
cross-element\/} from
$bT$ to $aT$ if $\gamma T=aT$ and $T\gamma=Tb$ (herein 
lies the significance
for us of matrix units). By including the standard 
one-dimensional infinite
shift among our generators, we can ensure that $A$ 
contains a cross-element
from any standard
primitive idempotent of $J_k$ to any one of $J_{k+1}$, 
also in $\roman
O((\log k)^2)$ products of generators.

Now suppose $\theta\:A\rightarrow G(s)$ is an algebra 
embedding for some $s<1$.
To complete the proof, we show, in several steps, that 
$r\leq s$. For each
$c\in \bold R^+$ let
$$W_s(c)=\{x\in B(F)\,|\,cn^s\text{ is a growth curve for 
}x\}\.$$
Notice that $G(s)$ is the union of these subspaces 
$W_s(c)$. Because $s<1$, the
powers of $W_s(c)$ grow ``polynomially''.
\rem{Step \rm1} For any $c\in \bold R^+$, there exists 
$d\in \bold R^+$ such
that for all $m\in \bold N$
$$(W_s(c))^m\subseteq W_s(dm^{1/1-s})\.$$
(We remark that for $s=0$, the $W_s(c)$ actually give a 
filtering of $G(0)$.
At the other extreme, for $s=1$ the powers of $W_s(c)$ grow
exponentially.)\endrem
\medskip

For each $k\in \bold N$, fix a set $F_k$ of $n_k$ 
equivalent orthogonal nonzero
idempotents of $J_k$ (for example, the standard 
primitives). From the above key 
property and Step 1, we obtain:
\rem{Step \rm 2} There exist nondecreasing $c_k\in \bold 
R^+$ for
$k=1,2,3\dotsc$ such that
\roster
\item \<$c_kn^s$ is a growth curve for elements in 
$\theta(F_k)$, for a
cross-element between each pair, and for a cross-element 
between any element of
$\theta(F_k)$ and one of $\theta(F_{k+1})$.
\item \<$c_k=\roman O((\log k)^{2/(1-s)})$.
\endroster
\endrem

Our strategy now is to chase the first nonzero rows of the 
images of the
idempotents in $F_k$ and obtain opposing constraints on 
their positions. On the
one hand, they have to be fairly close together to ensure 
that some
cross-element from $\theta(F_k)$ to $\theta(F_{k+1})$ can 
lie in $W_s(c_k)$. On
the other hand, having the $n_k$ images of the equivalent 
orthogonal
idempotents from $F_k$ all inside $W_s(c_k)$ forces their 
nonzero rows to
become increasingly scattered. A comparison of these 
opposing growths leads to:

\rem{Step \rm3} $n_k=\roman O(k^{s/1-s}c_k^{(1+s)/(1-s)})$.
\endrem

From (2) of Step 2, this yields $n_k=\roman 
O(k^{s/1-s}(\log k)^u)$ where
$u=
2(1+s)/(1-s)^2$. But by definition, $n_k$ has polynomial 
growth of $k^t$.
Hence $t\leq s/(1-s)$, which implies $r/(1-r)\leq s/(1-s)$ 
and so $r\leq s$, as
required.


\Refs
\rc
\ref\no1
\by K. R. Goodearl, P. Menal, and J. Moncasi
\paper Free and residually artinian regular rings
\jour J. Algebra
\vol 156
\yr 1993
\pages 407--432
\endref
\ref\no2
\by G. R. Krause and T. H. Lenagan
\book Growth of algebras and Gelfand-Kirillov dimension
\publ Pitman
\publaddr New York
\yr 1985
\endref
\ref\no3
\by K. C. O'Meara, C. I. Vinsonhaler, and W. J. Wickless
\paper Identity-preserving embeddings of countable rings 
into \RM2-generator
rings
\jour Rocky Mountain J. Math.
\vol 19
\yr 1989
\pages 1095--1105
\endref
\endRefs
\enddocument